\documentclass{amsart}
\usepackage{amssymb,amsmath}
\usepackage[mathscr]{euscript}
\input{epsf}

\newcounter{punct}

\def\punct{\refstepcounter{punct}{
\arabic{punct}.  }}

\newtheorem{theorem}{Theorem}

\newtheorem{lemma}[theorem]{Lemma}

          \def\sm{\smallskip}

\begin{document}

\newcommand{\supp}{\mathop {\mathrm {supp}}\nolimits}
\newcommand{\rk}{\mathop {\mathrm {rk}}\nolimits}
\newcommand{\Aut}{\mathop {\mathrm {Aut}}\nolimits}
\newcommand{\Ob}{\mathop {\mathrm {Ob}}\nolimits}
\newcommand{\Out}{\mathop {\mathrm {Out}}\nolimits}
\renewcommand{\Re}{\mathop {\mathrm {Re}}\nolimits}
\newcommand{\Inn}{\mathop {\mathrm {Inn}}\nolimits}
\newcommand{\Char}{\mathop {\mathrm {Char}}\nolimits}
\newcommand{\ch}{\cosh}
\newcommand{\sh}{\sinh}
\newcommand{\Sp}{\mathop {\mathrm {Sp}}\nolimits}
\newcommand{\SOS}{\mathop {\mathrm {SO^*}}\nolimits}
\newcommand{\Ams}{\mathop {\mathrm {Ams}}\nolimits}
\newcommand{\Gms}{\mathop {\mathrm {Gms}}\nolimits}

\newcommand{\edge}{\mathop {\mathrm {edge}}\nolimits}

\def\0{\mathbf 0}

\def\ov{\overline}
\def\wh{\widehat}
\def\wt{\widetilde}
\def\pol{\twoheadrightarrow}

\renewcommand{\rk}{\mathop {\mathrm {rk}}\nolimits}
\renewcommand{\Aut}{\mathop {\mathrm {Aut}}\nolimits}
\renewcommand{\Re}{\mathop {\mathrm {Re}}\nolimits}
\renewcommand{\Im}{\mathop {\mathrm {Im}}\nolimits}
\newcommand{\sgn}{\mathop {\mathrm {sgn}}\nolimits}

\def\bfa{\mathbf a}
\def\bfb{\mathbf b}
\def\bfc{\mathbf c}
\def\bfd{\mathbf d}
\def\bfe{\mathbf e}
\def\bff{\mathbf f}
\def\bfg{\mathbf g}
\def\bfh{\mathbf h}
\def\bfi{\mathbf i}
\def\bfj{\mathbf j}
\def\bfk{\mathbf k}
\def\bfl{\mathbf l}
\def\bfm{\mathbf m}
\def\bfn{\mathbf n}
\def\bfo{\mathbf o}
\def\bfp{\mathbf p}
\def\bfq{\mathbf q}
\def\bfr{\mathbf r}
\def\bfs{\mathbf s}
\def\bft{\mathbf t}
\def\bfu{\mathbf u}
\def\bfv{\mathbf v}
\def\bfw{\mathbf w}
\def\bfx{\mathbf x}
\def\bfy{\mathbf y}
\def\bfz{\mathbf z}

\def\bfA{\mathbf A}
\def\bfB{\mathbf B}
\def\bfC{\mathbf C}
\def\bfD{\mathbf D}
\def\bfE{\mathbf E}
\def\bfF{\mathbf F}
\def\bfG{\mathbf G}
\def\bfH{\mathbf H}
\def\bfI{\mathbf I}
\def\bfJ{\mathbf J}
\def\bfK{\mathbf K}
\def\bfL{\mathbf L}
\def\bfM{\mathbf M}
\def\bfN{\mathbf N}
\def\bfO{\mathbf O}
\def\bfP{\mathbf P}
\def\bfQ{\mathbf Q}
\def\bfR{\mathbf R}
\def\bfS{\mathbf S}
\def\bfT{\mathbf T}
\def\bfU{\mathbf U}
\def\bfV{\mathbf V}
\def\bfW{\mathbf W}
\def\bfX{\mathbf X}
\def\bfY{\mathbf Y}
\def\bfZ{\mathbf Z}

\def\frD{\mathfrak D}
\def\frQ{\mathfrak Q}
\def\frS{\mathfrak S}
\def\frT{\mathfrak T}
\def\frL{\mathfrak L}
\def\frG{\mathfrak G}
\def\frb{\mathfrak b}
\def\frg{\mathfrak g}
\def\frh{\mathfrak h}
\def\frf{\mathfrak f}
\def\frk{\mathfrak k}
\def\frl{\mathfrak l}
\def\frm{\mathfrak m}
\def\frn{\mathfrak n}
\def\fro{\mathfrak o}
\def\frp{\mathfrak p}
\def\frq{\mathfrak q}
\def\frr{\mathfrak r}
\def\frs{\mathfrak s}
\def\frt{\mathfrak t}
\def\fru{\mathfrak u}
\def\frv{\mathfrak v}
\def\frw{\mathfrak w}
\def\frx{\mathfrak x}
\def\fry{\mathfrak y}
\def\frz{\mathfrak z}

\def\bfw{\mathbf w}

\def\R {{\mathbb R }}
 \def\C {{\mathbb C }}
  \def\Z{{\mathbb Z}}
  \def\H{{\mathbb H}}
\def\K{{\mathbb K}}
\def\N{{\mathbb N}}
\def\Q{{\mathbb Q}}
\def\A{{\mathbb A}}

\def\T{\mathbb T}
\def\P{\mathbb P}
\def\SS{\mathbb S}

\def\G{\mathbb G}

\def\cD{\EuScript D}
\def\cL{\EuScript L}
\def\cK{\EuScript K}
\def\cM{\EuScript M}
\def\cN{\EuScript N}
\def\cP{\EuScript P}
\def\cT{\EuScript T}
\def\cQ{\EuScript Q}
\def\cR{\EuScript R}
\def\cW{\EuScript W}
\def\cY{\EuScript Y}
\def\cF{\EuScript F}
\def\cG{\EuScript G}
\def\cZ{\EuScript Z}
\def\cI{\EuScript I}
\def\cB{\EuScript B}
\def\cA{\EuScript A}
\def\cE{\EuScript E}
\def\cC{\EuScript C}
\def\cS{\EuScript S}

\def\bbA{\mathbb A}
\def\bbB{\mathbb B}
\def\bbD{\mathbb D}
\def\bbE{\mathbb E}
\def\bbF{\mathbb F}
\def\bbG{\mathbb G}
\def\bbI{\mathbb I}
\def\bbJ{\mathbb J}
\def\bbL{\mathbb L}
\def\bbM{\mathbb M}
\def\bbN{\mathbb N}
\def\bbO{\mathbb O}
\def\bbP{\mathbb P}
\def\bbQ{\mathbb Q}
\def\bbS{\mathbb S}
\def\bbT{\mathbb T}
\def\bbU{\mathbb U}
\def\bbV{\mathbb V}
\def\bbW{\mathbb W}
\def\bbX{\mathbb X}
\def\bbY{\mathbb Y}

\def\kappa{\varkappa}
\def\epsilon{\varepsilon}
\def\phi{\varphi}
\def\le{\leqslant}
\def\ge{\geqslant}

\def\B{\mathrm B}

\def\la{\langle}
\def\ra{\rangle}
\def\tri{\triangleright}

\def\lambdA{{\boldsymbol{\lambda}}}
\def\alphA{{\boldsymbol{\alpha}}}
\def\betA{{\boldsymbol{\beta}}}
\def\mU{{\boldsymbol{\mu}}}

\def\const{\mathrm{const}}
\def\rem{\mathrm{rem}}
\def\even{\mathrm{even}}
\def\SO{\mathrm{SO}}
\def\SL{\mathrm{SL}}
\def\PSL{\mathrm{PSL}}
\def\cont{\mathrm{cont}}
\def\O{\mathrm{O}}
\def\Th{\mathrm{Th}}

\def\U{\operatorname{U}}
\def\GL{\operatorname{GL}}
\def\Mat{\operatorname{Mat}}
\def\End{\operatorname{End}}
\def\Mor{\operatorname{Mor}}
\def\Aut{\operatorname{Aut}}
\def\inv{\operatorname{inv}}
\def\red{\operatorname{red}}
\def\Ind{\operatorname{Ind}}
\def\dom{\operatorname{dom}}
\def\im{\operatorname{im}}
\def\md{\operatorname{mod\,}}
\def\indef{\operatorname{indef}}
\def\Gr{\operatorname{Gr}}
\def\Pol{\operatorname{Pol}}
\def\Cut{\operatorname{Cut}}
\def\Add{\operatorname{Add}}
\def\ord{\operatorname{ord}}
\def\Replace{\operatorname{Replace}}
\def\Tr{\operatorname{Tr}}
\def\Homeo{\operatorname{Homeo}}
\def\Sep{\operatorname{Sep}}

\def\arr{\rightrightarrows}
\def\bs{\backslash}

\def\cH{\EuScript{H}}
\def\cO{\EuScript{O}}
\def\cQ{\EuScript{Q}}
\def\cL{\EuScript{L}}
\def\cX{\EuScript{X}}

\def\Di{\Diamond}
\def\di{\diamond}

\def\fin{\mathrm{fin}}
\def\ThetA{\boldsymbol {\Theta}}

\def\krest{{\text{\tiny{\rm\texttimes}}}}
\def\vph{{\vphantom{\bigr|}}}

\begin{center}
\large\bf
 On the group of homeomorphisms of the Basilica
 
 \medskip

\sc Yury A.Neretin%
	\footnote{Supported by the grant FWF,  P31591.}
 
\end{center}

\medskip

{\small We show that the group  of all homeomorphisms of the Basilica fractal
coincides with a group of transformations of a certain non-locally finite ribbon tree.
Also, we show that Basilica Thompson group defined by Belk and Forrest is dense in the group
of all orientation preserving homeomorphisms of the Basilica.

}

\medskip

{\bf \punct The Basilica.} The Basilica $\cB$ is a compact topological space obtained as a quotient of the circle $S=\R/2\Z$ under the equivalence relation
\begin{equation}
\tfrac{3k+1}{3^n}\sim \tfrac{3k+2}{3^n}, \qquad \text{where $k$, $n$ are nonnegative integers,}
\label{eq:eq}
\end{equation}
In particular $1\sim 2=0$.
In other words, we write  numbers in ternary system and identify $a_0.a_1\dots a_{n-1}1\sim a_0.a_1\dots a_{n-1}2$.
The reader can  find nice pictures of the Basilica in papers \cite{LM}, \cite{BF1}--\cite{BF2} or their preprint versions in {\tt arXiv},
see also  Fig.~1.

We say that a {\it separating point} of $\cB$ is a point corresponding a nontrivial (two-point) equivalence class.
Denote by $\Sep(\cB)$ the set of separating points. If $x\in\Sep(\cB)$, then $\cB\setminus x$ consists of two components.
Otherwise, $\cB\setminus x$ is connected.

It is convenient
to think that our circle $S$ is a boundary of a disk $D$, and equivalent points are connected by  chords $[a\text{---}b]$ ({\it leafs}), clearly they are disjoint.
Denote the set of all leafs by $L(\cB)$.
Leafs divide the disc into domains, say {\it Cantorgons} (cf. \cite{Det}), having a countable number of sides, an intersection of a Cantorgon with $S$ 
is a totally disconnected (Cantor) set. The set of sides of a Cantorgon has a canonical cyclic order; as a cyclic ordered set this set
 is equivalent to $\Q/\Z$ (or to a countable dense subset of a circle, all such subsets are equivalent).

\sm

{\sc Remark.}
The set of all leafs can be obtained as follows. Consider a circle of length 2. Its diameter $[0\text{---}1]$ splits a circle into two arcs.
We divide each arc by a pair of points with ratios $1:1:1$ and connect  elements of each pair by a chord. We get 6 ends of chords and 6 arcs and divide arcs with ratios
$1:1:1$, etc.
\hfill $\boxtimes$

\sm

Consider a Cantorgon $C$ and contract each side to a point. Then we get a topological disc, its boundary
$\eth(C)$ 
is a topological circle equipped with a canonical orientation induced from $S$. Each circle $\eth(C_j)$ is homeomorphically  embedded to $\cB$,
we call such image  a {\it component} of the Basilica.
Two components are disjoint or contain one joint point, which with necessity is a separating point. The union of all components is dense in the Basilica.

Clearly, homeomorphisms of $\cB$ send components to components.
We say that a homeomorphism of $\cB$ {\it preserves orientation} if it is an orientation preserving map on each component.  
Denote by $\Homeo(\cB)$ the group of all homeomorphisms of $\cB$, by $\Homeo_+(\cB)$ the subgroup of orientation preserving homeomorphisms.

\sm

{\bf\punct The tree of the Basilica.}
Consider a graph $T$, whose vertices  are enumerated by Cantorgons of $\cB$ and edges are enumerated by leafs; an edge  contains  a vertex
 if the corresponding leaf is a side of the corresponding  Cantorgon. Clearly, $T$ is a tree and valences of all vertices are countable.
 Denote by $\edge(v)$ the set of edges incident to a vertex $v$. As we have seen above, each set $\edge(v)$ has a canonical cyclic order.
 Denote the tree equipped with this structure by $T(\cB)$. 
 
 Consider an automorphism $g$ of the tree $T$ (i.e., a bijection of the set of vertices to itself sending adjacent vertices to adjacent vertices).
 We say that $g$ is {\it orientation preserving} if for each vertex $v$ the map $g$ induces a cyclic order preserving map
 $\edge(v)\to \edge(gv)$. We say that $g$ {\it respects orientation} if each induced map  $\edge(v)\to \edge(gv)$ preserves or inverses the cyclic order.
 We denote by $\Aut_+(T(\cB))$  (resp., $\Aut(T(\cB))$) the groups  of all orientation preserving (resp., orientation respecting)
 automorphisms of $T(\cB)$ (this construction is inside the formalism of Burger--Mozes groups, see \cite{BM}, Subsect. 3.2).
 
 \begin{figure}
  $$\epsfbox{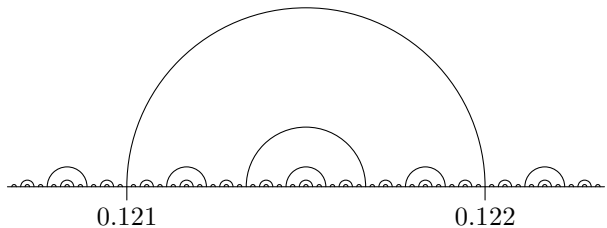}$$
  \caption{We draw a circle $\R/2\Z$ as a segment  of the horizontal line and connect equivalent points by arcs.}
 \end{figure}

 We define topology on these group by the assumption: stabilizers of finite subtrees
 are open subgroups. It is easy to see that we get  topological groups, and, moreover,  Polish groups (a group is called {\it Polish} 
 if 
 it is homeomorphic to a complete metric space%
 \footnote{Denote by $E$ the set of edges of $T$, by $V$ the set of vertices, by $S(E)$, $S(T)$ the groups of all permutations
 of these sets equipped with the unique separable topology on infinite symmetric group, see, e.g., \cite{Gao}. Let $A\subset V\times E$ consists of pairs $(v,e)$ such that $v\in e$. Let $B\subset V\times E\times E\times E$
 consists of $(v,e_1,e_2,e_3)$ such that $v\in e_j$ and $e_1$, $e_2$, $e_3$ are arranged clockwise. Then $\Aut_+(T(\cB))$
 is the subgroup in a Polish group $S(E)\times S(V)$ consisting of elements sending $A$ to $A$ and $B$ to $B$. It is easy to see that this condition is closed,
 and therefore $\Aut_+(T(\cB))$ is a closed subgroup in a Polish group, so it is Polish.}).

\sm

{\bf\punct The statement of the paper.}
Since a homeomorphism of the Basilica sends components to components and separating points to separating points,
we have a canonical homomorphism 
 $$
 \Pi:\Homeo(\cB)\to \Aut(T(\cB)),
 $$
 sending $\Homeo_+(\cB)\to \Aut_+(T(\cB))$.  Since separating points are dense, our homomorphism is injective.
 
 
 \begin{theorem}
 \label{th}
 {\rm a)} $\Pi$ establishes an isomorphism
  $ \Homeo(\cB)\simeq \Aut(T(\cB))$.
  
  \sm
 {\rm b)} $\Pi$ establishes an isomorphism $ \Homeo_+(\cB)\simeq \Aut_+(T(\cB))$.
  
  \sm
  
  {\rm c)} The Basilica Thompson  group $\Th(\cB)$ is dense in $\Aut_+(T(\cB))$.
 \end{theorem}
 
 
 \sm
 
 {\bf \punct The  Basilica Thompson group $\Th(\cB)$.} 
 We say that a {\it diagram} is a finite collection $\Delta\subset L(\cB)$ of leafs such that 
 
 \sm
 
 1) for any $\eta\in L(\cB)\setminus \Delta$
 all elements of $\Delta$ lie on one side of $\eta$.

 \sm
 
 2) $\Delta$ contains the leaf $[0\text{---}1]$.
 
 \sm
 
 Leafs of $\Delta$ divide the disk $D$ into domains, say {\it cells}. 
 Our condition means 
 that for each cell all adjacent leafs are sides of one Cantorgon.  We call ends of leafs {\it vertices of a diagram},
 and arcs of $S$ between vertices {\it intervals of a diagram}.

 We say that a diagram is {\it perfect}
 if it   can be obtained by the following procedure.
 We start with the diameter $[1-2]$. It splits $S$ into two intervals. We divide one of intervals with ratios $1:1:1$ by points $a$, $b$ and draw the corresponding chord $[a-b]$.
 Now a circle is split into 4 intervals. We choose one interval, divide it with ratios $1:1:1$ and add the corresponding chord. Etc.
 By the construction, intervals of a perfect diagram have forms
  \begin{equation}
   \bigl[\tfrac{3\alpha}{3^a},\tfrac{3\alpha+1}{3^a}\bigr],\qquad   \bigl[\tfrac{3\beta+1}{3^b},\tfrac{3\beta+2}{3^b}\bigr],
   \qquad  \bigl[\tfrac{3\gamma+2}{3^c},\tfrac{3\gamma}{3^c}\bigr].
  \end{equation}
  
We say that two perfect diagrams $\Delta_1$, $\Delta_2$ are isomorphic if there is an orientation preserving bijective continuous map identifying 
$\Delta_1$ and $\Delta_2$. An isomorphism determines a bijective map on the set of vertices. We extend this map to each interval of $\Delta_1$
in a linear way (all slopes have a form $3^k$). According Belk and Forrest  \cite{BF1}, such maps form a group---the {\it Basilica Thompson group} $\Th(\cB)$, see \cite{BF1} (for instance, it is generated by 4 explicit elements, Theorem 7.1), for further properties, see \cite{LM}, \cite{WZ}, see also \cite{BF2} on  extensions of this construction. On the other hand, such a map respects the equivalence relation
\eqref{eq:eq}
and therefore determines a homeomorphism of the Basilica.

\sm

{\bf\punct Correspondence with the initial model of the Basilica Thompson group.} The Basilica naturally arises as an object of
complex dynamics (as a Julia set, see, e.g., \cite{LM}), a natural equivalence relation on a circle is defined in the following way. We consider a circle $S^\diamond=\R/\Z$
of length 1 and the corresponding disc. We start with chords $[\frac 16\text{---}(-\frac16)]$, $[\frac 13\text{---}\frac23]$ and apply the procedure
described in Remark in Sect.1 changing ratios $1:1:1$ to $1:2:1$. In this way, 
 we get the same combinatorial picture and the same
topological object. Again, we consider perfect diagrams, their isomorphisms,
and linear maps on intervals (now all slopes have form $2^k$).

Let us describe the correspondence more precisely. First,
we pass to the circle $S^\circ=\R/3\Z$ with chords $[0\text{---}1]$,
$[3/2\text{---}5/2]$ (we rotate the circle $S^\diamond$ and dilatate it).
The homeomorphism $\Theta: S\to S^\circ$ is defined in the following way. Let $a$ range in the set $\{0,1,2\}$. We define
functions $\phi(a)$, $\psi(a)$
by
$$
\phi(0)=\tfrac14,\, \phi(1)=\tfrac 12,\, \phi(2)=\tfrac14,\qquad \psi(0)=0,\,\psi(1)=\tfrac14,\,\psi(2)=\tfrac 34. 
$$
Then for a ternary number $A=a_0.a_1a_2\dots$ we set
$$
\Theta(A)=\frac{\psi(a_0)}{\phi(0)}+ \sum_{n=1}^\infty \Bigl(\prod_{p=0}^{n-1} \phi(a_p)\cdot \psi(a_n)\Bigr).
$$
This homeomorphism identifies the equivalence relations and identify the group
$\Th(\cB)$ defined above with the group from \cite{BF1}.

\sm

{\bf\punct Identification of $\Homeo_+(\cB)$ and $\Aut_+ T(\cB)$.}

\begin{lemma}
\label{l:extension}
 Let $\Delta$ be a diagram. Consider all cells $M_j$ of $\Delta$ adjacent to $>1$ leafs
 and the corresponding Cantorgons $C_j$. Then there is a finite collection of leafs
 $l_1,l_2,\dots$, such that each $l_q$ is a side of some Cantorgon $C_p$ and the diagram
 $\Delta\cup \{l_1,l_2,\dots\}$ is perfect.
\end{lemma}

{\sc Proof.} Let a leaf $[a_0.a_1\dots a_n 1\text{---}a_0.a_1\dots a_n 2]$
be contained in a perfect diagram $\Sigma$. It can appear as a result of cutting of the arc
$[a_0.a_1\dots a_n,a_0.a_1\dots a_n+3^{-n}]$ with ratios $1:1:1$ and therefore ends of this arc are vertices
of $\Sigma$, i.e., ends of arcs of $\Sigma$. Minimal ternary notation of these vertices has $\le n$ places, 
therefore the corresponding leafs are longer than our initial leaf.

Consider the diagram $\Delta$. For each cell $M_j$ we find an adjacent leaf $m_j$ of minimal length and add to $\Delta$ all
sides of the Cantorgon $C_j$ whose lengths are larger  than the length of $m_j$. We get a perfect diagram. 
\hfill $\square$

\begin{lemma}
\label{l:3}
Let $r\in\Aut_+(T(\cB))$,
let $P$ be a finite subtree in $T(\cB)$.
 Then there exists  $g\in \Th(\cB)$ such that $ g\Bigr|_{P}=\Pi\circ r\Bigr|_P$.
\end{lemma}

{\sc Proof.} Without loss of generality we can assume that both $P$, $rP$ contain
the edge corresponding the leaf $[0\text{---}1]$.

For each vertex of $P$ the set of adjacent vertices is equipped with a cyclic order.
So $P$ has a structure of a planar tree; $rP$ is an isomorphic planar tree. 
For each edge of $P$ we draw the corresponding leaf and get a diagram, say $\Delta$. 
Let $\Sigma$ be the similar diagram for $rP$. These diagrams are combinatorially  isomorphic.
Next, we extend them to perfect diagrams $\wt \Delta$, $\wt\Sigma$ as in Lemma
\ref{l:extension}. This means that we added additional leafs between each pair of neighbouring  leafs of each cell
of $\Delta$ (resp. $\Sigma$). Generally, diagrams $\wt \Delta$, $\wt \Sigma$ are non-isomorphic since number of additional leafs
for corresponding pairs are different. However,   adding more leafs to $\wt \Delta$ and $\wt \Sigma$
we can get isomorphic perfect diagrams. It remains to   take the corresponding $g\in\Th(\cB)$.
\hfill $\square$

\sm

{\sc Proof of statements b), c) of Theorem \ref{th}.}
Fix $r\in \Aut_+(T(\cB))$ 
Consider a sequence of finite subtrees $P_j\subset T(\cB)$ such that
$\cup P_j= T(\cB)$. For each $P_j$, consider an element $g_j\in \Aut_+(T(\cB))$  
such that $ g_j\Bigr|_{P_j}=\Pi\circ r\Bigr|_{P_j}$.

Then for each $x_\alpha\in\Sep(\cB)$ the sequence $g_j (x_\alpha)$ is eventually constant, denote
its limit by $\rho(x_\alpha)$.
For each $m$ denote by $\Omega_m$ the set of $x_\alpha\in \Sep(\cB)$ such that
$g_jx_\alpha=\rho(x_\alpha)$  for $j\ge m$. Clearly, we have $\cup \Omega_m=\Sep(\cB)=\cup\rho(\Omega_m)$.
For this reason, the map $\rho:\Sep(\cB)\to\Sep(\cB)$ is a bijection. 

Next, sets $\cB\setminus \Omega_m$, $\cB\setminus \rho(\Omega_m)$ split into finite collection of connected
pieces and diameters of these pieces tend to 0 as $m\to\infty$. Elements $g_j$ for $j\ge m$ send  
pieces to pieces. This implies uniform convergence of homeomorphisms $g_j$. We define $\rho$ as 
limit of $g_j$, the corresponding transformation of the tree $T(\cB)$ is $r$.

Thus, we established that $\Th(\cB)$ is dense in $\Homeo_+(\cB)$ and  the homomorphism
$\Pi:\Homeo_+(\cB)\to \Aut_+(T(\cB))$ is surjective. Continuity of $\Pi$ is obvious. Banach's bounded inverse theorem 
remains true for Polish groups, see e.g., \cite{Gao}, Corollary 2.3.4. For this reason, $\Pi$ is an isomorphism
of Polish groups.
\hfill 

\sm

{\bf\punct Identification of $\Homeo(\cB)$ and $\Aut T(\cB)$.}
For any  leaf $[A\text{---}A']=[a_0.a_1\dots a_n1\text{---}a_0.a_1\dots a_n2]$ consider 
the following  homeomorphism $\gamma_{A,A'}$ of $\cB$: 
it permutes $A$ and $A'$, it is linear on the interval $[A,A']$, and
is identical outside this interval.

We can repeat the previous considerations for the countable group $\Gamma$ generated by $\Th(\cB)$ and $\gamma_{A,A'}$. Indeed, denote by $e_{A,A'}$ the corresponding edge of $T(\cB)$. Cutting it in a midpoint we splits the tree into two branches. The transformation
$\gamma_{A,A'}$ is identical on one of branches and is a map changing cyclic ordering in each vertex of another branch. Now we can get the analogue
of Lemma \ref{l:3}. Let $q\in \Aut(T(\cB))$, $P$ be a finite subtree. We can
choose a collection $(A_1, A'_1)$, \dots, $(A_k, A'_k)$ such that
$$
r=\gamma_{A_k, A'_k}\dots \gamma_{A_1, A'_1} q\in \Aut_+(T(\cB))
$$
We apply Lemma \ref{l:3} and get an element $g\in\Gamma$ such that
$g=\Pi\circ q\Bigr|_P$. Further arguments are identical.

\sm

{\bf \punct Final remarks.} a) I have not seen the groups $\Aut_+(T(\cB))$ and 
$\Aut(T(\cB))$ in literature.  However, these objects are not as exotic as it may seem at first glance. In 1982 Olshanski
\cite{Olsh} classified  unitary representations of
the group 
 of all automorphisms of a non-locally finite
Bruhat--Tits tree, in a certain sense this example is toy, but  this work was an important test site for representation theory of infinite-dimensional groups. The group
of spheromorphisms of this tree is a nontrivial object for representation theory
\cite{Ner}.

The group of order preserving bijections of the set $\Q$ is a standard example  of  Polish groups,
see, e.g., \cite{Pes}, \cite{Tsa}, \cite{JTs}. It also admits explicit classification of unitary representations \cite{Tsa} (but in this case the list of representations is  poor). On the other hand, dense cyclic orderings participate in the construction of 'Chinese restaurant process' by Pitman, see, e.g.,
\cite{Pit},
this process is a starting point for harmonic analysis related to 
  infinite symmetric groups.

\sm

b) Belk and Forrest \cite{BF2} proposed a general construction of discrete groups of Thompson type acting on certain fractal
compact spaces. Apparently, there are  statements similar to Theorem \ref{th} for some of these spaces. Certainly,
this is so for the Doaudy rabbit known in complex dynamics. Another interesting possible example is the Airplane fractal defined in \cite{BF2}, see also \cite{Tar}.
Apparently, in this case we must consider the  Bruhat--Tits $\R$-tree
(see, e.~g., \cite{Mor}) of the field of formal Laurent series 
$\Q[[t,t^{1,2}, t^{1/4}, \dots]]$ equipped with a ribbon structure.
Apparently, the group of homeomorphisms of the Airplane is isomorphic to the group of {\it homeomorphisms}
of this $\R$-tree
 respecting the cyclic ordering structure of germs of edges in vertices.
 
 \sm 
 
 c) It is well-known that unitary representation of discrete groups are a dangerous topic. However, classification \cite{Olsh}
 of unitary representations of a non-locally finite Bruhat--Tits tree
   contains a counterpart of the complementary series. Restricting it to $\Th(\cB)$ we get unitary representations that are at least unusual and are spherical with respect to the stabilizer of a Cantorgon (this stabilizer is described in
 \cite{BF1}, Sect.~6).

\sm

I am grateful to Matteo Tarocchi for a discussion of $\Th(\cB)$.

\tt Fakult\"at f\"ur Mathematik, Universit\"at Wien;

Institute for Information Transmission Problems;

Moscow State University

e-mail: yurii.neretin(frog)univie.ac.at

URL: mat.univie.ac.at/$\sim$neretin

\end{document}